\documentstyle[12pt]{article}
\normalbaselines

\begin{document}
\begin{flushright}

{\raggedleft

math.HO/9904077 v2\\[1cm]}

\end{flushright}

\renewcommand{\thefootnote}{\fnsymbol{footnote}}
\begin{center}
{\LARGE\baselineskip0.9cm 
Angles in Complex Vector Spaces\\[1.5cm]}

{\large K. Scharnhorst\footnote[2]{E-mail:
{\tt scharnh@physik.hu-berlin.de}}
}\\[0.3cm]

{\small Humboldt-Universit\"at zu Berlin

Institut f\"ur Physik

Invalidenstr.\ 110

D-10115 Berlin

Federal Republic of Germany}\\[1.5cm]

\begin {abstract}
The article reviews some of the (fairly scattered) 
information available in the
mathematical literature on the subject of angles in complex vector spaces.
The following angles and their relations are considered: 
Euclidean, complex, and Hermitian angles, (Kasner's) pseudo-angle, 
the K\"ahler angle (synonyms for the latter used in the literature are:
angle of inclination, characteristic deviation, holomorphic
deviation, holomorphy angle, Wirtinger angle, slant angle).
\end{abstract}

\end{center}

\renewcommand{\thefootnote}{\arabic{footnote}}

\thispagestyle{empty}

\newpage

\noindent
{\bf 1.\ Introduction.}\ The angle between two vectors in a 
real vector space is a concept often already introduced to
students at the school level. Complex vector spaces feature 
prominently in most linear algebra courses at the undergraduate level
and they can be found in many branches of mathematics, the natural
sciences, and engineering. Therefore, it is surprising that very little 
guidance is available in the mathematical literature on 
angles in complex vector spaces. Even advanced books on linear
algebra or higher geometry hardly mention the subject. In order
to find some information on it (which is widely scattered, however) 
one has to resort almost exclusively to the journal literature. 
No review seems to be 
available, closest to this come only sections in the recently 
published monographs by Rosenfeld 
(\cite{roze2}, chap.\ III, \S 3.3, sect.\ 3.3.6, pp.\ 182/183)
and, more detailed, by Goldman (\cite{gold}, chap.\ 2, sect.\ 2.2.2,
pp.\ 36-39) which, however, are supplemented by no or only
very few references, respectively, concerning the subject (the latter 
monograph appeared in print only when the first version of the present
paper had been put into circulation on the Los Alamos Math Archive). 
In the present article, which
grew from the working needs of a theoretical physicist, we undertake
to fill this gap for a rather general audience.\\

To begin with, let us consider the problem one faces in 
introducing the concept of an angle in complex vector spaces.
In any (finite-dimensional) real (Euclidean) vector space $V_{\bf R}$ 
($\simeq {\bf R}_m$, $m \in {\bf N}$, $m\ge 2$) equipped with the
scalar product $(A,B)_{\bf R} = \sum_{k=1}^m A_k B_k$ 
for any pair of vectors $A, B \in V_{\bf R}$ 
one can define an (real) angle $\Theta(A,B)$, 
$0\le\Theta\le \pi$,
between these two vectors 
by means of the standard formula 
($\vert A\vert = \sqrt{ (A,A)_{\bf R}}$)
\begin{eqnarray}
\label{ZA1}
\cos \Theta(A,B)&=& 
\frac{\ (A,B)_{\bf R}}{\vert A\vert\ \vert B\vert}\ \ .
\end{eqnarray}
The introduction of an angle between two vectors $a, b$ of
a (finite-dimensional) complex (Hermitian, unitary) 
vector space $V_{\bf C}$ 
($\simeq {\bf C}_n$, $n \in {\bf N}$, $n\ge 2$) is ambiguous
and can be performed
\begin{itemize}
\item[I.]
either directly in the complex
vector space $V_{\bf C}$ by relying on the Hermitian product
$(a,b)_{\bf C} = \sum_{k=1}^n \bar{a}_k b_k$ defined in it
for any pair of vectors $a, b \in V_{\bf C}$  
($\bar{a}_k$ denotes the complex conjugate of $a_k\in {\bf C}$), or

\item[II.]
by relying on the real vector space $V_{\bf R}$ 
($\simeq {\bf R}_{2n}$) isometric to $V_{\bf C}$.
\end{itemize}
Both approaches which are not completely independent 
are equally justified. Consequently, one has
to study their relation and this is what this article mainly
is concerned with.\\

\noindent
{\bf 2.\ The Euclidean, complex, Hermitian, and pseudo-angles.}\ 
The (real-valued) {\it Euclidean angle} between two vectors 
$a, b \in V_{\bf C}$ related to item II.\ 
and denoted by $\Theta (a,b)$, $0\le\Theta\le \pi$,
is defined by the formula
\begin{eqnarray}
\label{ZA2}
\cos \Theta (a,b)&=& \cos\Theta (A,B)\ =\
\frac{\ (A,B)_{\bf R}}{\vert A\vert\ \vert B\vert}\ \ ,
\end{eqnarray}
where we choose to determine the components of the 
vectors $A$, $B$ in the real vector space $V_{\bf R}$ 
by means of the relations
$A_{2k-1} = {\rm Re}\ a_k$ and $A_{2k} = {\rm Im}\ a_k$, 
$k = 1,\ldots,n$ \footnote{Throughout the article we denote
vectors in the real vector space $V_{\bf R}$ by capital letters
and vectors in the complex vector space $V_{\bf C}$ by small letters.}.
On the other hand, related to item I., a complex (-valued) angle 
$\Theta_{\rm c}(a,b)$ can be introduced by means of the
relation ($\vert a\vert = \sqrt{ (a,a)_{\bf C}} = \vert A\vert $)
\begin{eqnarray}
\label{ZA3}
\cos \Theta_{\rm c}(a,b)&=& 
\frac{\ (a,b)_{\bf C}}{\vert a\vert\ \vert b\vert}\ \ .
\end{eqnarray}
Here, we follow the definition applied in the majority of the 
literature on the subject (cf., e.\ g., \cite{roze6}, chap.\ VI,
\S 2, eq.\ (6.65), p.\ 574; 
incidentally, note that also other
views can be found, cf.\ \cite{frod}). Furthermore, one can 
write eq.\ (\ref{ZA3}) as
\begin{eqnarray}
\label{ZA4}
\cos \Theta_{\rm c}(a,b)&=& 
\rho\ {\rm e}^{\displaystyle i\varphi}\ \ , 
\end{eqnarray}
where ($\rho\le 1$)\footnote{This follows from the Cauchy
inequality $(a,b)_{\bf C} (b,a)_{\bf C}\le 
(a,a)_{\bf C} (b,b)_{\bf C}$ which can be derived from
$0\le (a-\sigma b,a-\sigma b)_{\bf C}$, 
$\sigma = (b,a)_{\bf C}/\vert b\vert ^2$ 
(\cite{roze6}, chap.\ VI, \S 2, eq.\ (6.67), p.\ 574).}
\begin{eqnarray}
\label{ZA5}
\rho &=& \cos\Theta_{\rm H}(a,b)\ =\ \vert \cos \Theta_{\rm c}(a,b)\vert\ \ .
\end{eqnarray}
$\Theta_{\rm H}(a,b)$, $0\le \Theta_{\rm H}\le \frac{\pi}{2}$,  
is called the {\it Hermitian angle}\footnote{This term (not the 
concept itself, however) has apparently been used first in
\cite{rizz4}, \S I., sect.\ 4, p.\ 96, also see \cite{brun6}, \cite{brun7}
(historically, there have been minor differences in the definition used
by some Italian authors, see, e.g., footnote 3 in \cite{brun6}, 
p.\ 395). For some early work on angles in complex vector spaces see 
\cite{fubi}-\cite{blas} (for further references see the reprint volume
cited in \cite{blas}).} between the vectors $a, b \in V_{\bf C}$ while 
$\varphi = \varphi(a,b)$, $-\pi\le\varphi\le\pi$ is called their
(Kasner's) {\it pseudo-angle} \cite{kasn1}-\cite{deci2}\footnote{
In a somewhat different context, the term 
{\it pseudo-angle} has already earlier been used by 
Giraud \cite{gira}, pp.\ 68/69, who calls this
quantity the {\it second 
pseudo-angle} while he denotes the Hermitian angle by the term
{\it first pseudo-angle} (see
\cite{gold}, chap.\ 2, sect.\ 2.2.2, p.\ 38).}. 
The pseudo-angle
$\varphi$ is of relevance in the context of pseudoconformal transformations. 
However, if one is interested just in the angle between two lines 
in the complex vector space $V_{\bf C}$ which are given by the 
vectors $a, b \in V_{\bf C}$ the concept of the pseudo-angle
can be disregarded. As any line in $V_{\bf C}$ defined by the 
vector $a$ can also be be given in terms of any other vector 
$a^\prime = z a$ ($z\in {\bf C}$,
$z\neq 0$) the pseudo-angle $\varphi$ does not have any 
meaning in this context. In other words, it can be disregarded in complex
projective spaces ${\bf CP}^n$. The Fubini-Study metric employed in
such spaces is just given by the cosine of the Hermitian angle (\ref{ZA5}).
The Hermitian angle can be understood 
geometrically as follows (\cite{shir1}, chap.\ IV, \S 32, pp.\ 405/406).
As in any real vector space the cosine of the (Hermitian) angle 
between two vectors $a,\ b\in V_{\bf C}$ can be defined to be the 
ratio between the length of the orthogonal projection
(with respect to the Hermitian product) of, say, the 
vector $a$ onto the vector $b$ to the length of the vector $a$ itself
(this projection vector is equal to 
$\sigma b$ where $\sigma = (b,a)_{\bf C}/\vert b\vert^2$)\footnote{Note, 
that there is a misprint at the given location in \cite{shir1}, p.\ 406,
9th line from above:  `perpendikulyara' has to be replaced 
by the term `vektora $x$'.}. From this definition 
immediately follows eq.\ (\ref{ZA5}).\\

\noindent
{\bf 3.\ The K\"ahler angle.}\ 
In order to proceed further let us introduce the almost complex structure
$J$, $J^2 = -1$, which acts as an operator
in the real vector space $V_{\bf R}$ isometric to $V_{\bf C}$. 
In our coordinates the almost complex structure $J$ performs
the following transformations:
$A_{2k-1}\longrightarrow A_{2k}$, $A_{2k}\longrightarrow -A_{2k-1}$,
$k = 1,\ldots,n$. This is equivalent to the transformation
$a\longrightarrow i a$ in $V_{\bf C}$. A subspace ${\cal P}$ of $V_{\bf R}$
is called {\it holomorphic}, if it holds ${\cal P} = J {\cal P}$. It is 
called {\it antiholomorphic\footnote{Sometimes, this term is also
used in a different sense.} (totally real, with a real 
Hermitian product)}, if it holds ${\cal P}\perp J{\cal P}$. Following the 
convention applied in a large fraction of the literature we 
introduce the notation $\tilde{A} = JA$, $A\in V_{\bf R}$
(note, $(\tilde{A},B)_{\bf R} = - (A,\tilde{B})_{\bf R}$
for any two vectors $A,\ B\in V_{\bf R}$).
By writing 
\begin{eqnarray}
\label{ZA6}
\cos\Theta_{\rm K}(a,b)\ \sin\Theta(a,b)
&=& \cos\Theta_{\rm K}(A,B)\ \sin\Theta(A,B)\ =\ 
\frac{\ (\tilde{A},B)_{\bf R}}{\vert A\vert\ \vert B\vert}\ \ \ \
\end{eqnarray}
one can now introduce a further
angle $\Theta_{\rm K}(a,b) = \Theta_{\rm K}(A,B)$,
$0\le \Theta_{\rm K}\le \pi$, which
is called the 
{\it K\"ahler angle}\footnote{
In using this term we follow the majority of the more recent 
mathematical literature (mostly, follow-up articles to \cite{cher}). 
The fundamental 2-form $\Phi(A,B) = (\tilde{A},B)_{\bf R}$ 
(\cite{koba}, chap.\ IX, \S 4, p.\ 147, also see \cite{yano}, \S4, p.\ 182,
eq.\ (4.15) and \cite{koba}, chap.\ IX, \S 7, p.\ 167
for the definition of the corresponding (K\"ahler) angle)
sometimes is referred to as the {\it K\"ahler function (K\"ahler form)}. 
Other terms which are or have been used for the K\"ahler angle are: 
{\it angle of inclination} \cite{shir2}, \S 3, p.\ 369,
\cite{roze5}, sect.\ 3, p.\ 378, \cite{roze3}, \cite{roze4},
{\it characteristic deviation} \cite{rizz1},
\cite{rizz4}, \S I., sect.\ 4, p.\ 96,
{\it holomorphic deviation} \cite{rizz2},
{\it holomorphy angle} \cite{roze2}, chap.\ III, \S 3.3, sect.\
3.3.6, pp.\ 182/183,
\cite{gold}, chap.\ 2, sect.\ 2.2.2, p.\ 36,
{\it Wirtinger angle}, {\it slant angle} \cite{chen1}-\cite{chen3}.} 
between the vectors $a,\ b\in V_{\bf C}$,
or the vectors $A,\ B\in V_{\bf R}$, respectively\footnote{If 
one wants to disregard the orientation of the 2-planes
in $V_{\bf R}$ one can take on the r.h.s.\ of eq.\ (\ref{ZA6})
the absolute value restricting the K\"ahler angle to the interval
$[0,\frac{\pi}{2}]$.}.
It is an intrinsic property of the (oriented) 2-plane in $V_{\bf R}$
defined by the (non-parallel) vectors $A$, $B$ with respect to the 
almost complex structure $J$ and, therefore, does not depend on 
the choice of the vectors in this plane  
(\cite{roze5}, sect.\ 3, p.\ 377, \cite{rizz1},
\cite{roze3}, sect.\ 1, p.\ 65 (p.\ 83 of the English
translation), \cite{roze2}, chap.\ III, \S 3.3, sect.\
3.3.6, p.\ 182)\footnote{For a 2-plane in $V_{\bf R}$ given by two 
unit vectors $a$, $b$ orthogonal in the Euclidean sense
($\Theta (a,b) = \frac{\pi}{2}$, $(a,b)_{\bf C} = -(b,a)_{\bf C}$) 
one can convince oneself that
$(a,b)_{\bf C}$ is invariant under the transformation
$a^\prime = \cos\phi\ a+ \sin\phi\ b$, 
$b^\prime = -\sin\phi\ a+ \cos\phi\ b$, $\phi\in {\bf R}$.}. 
The K\"ahler angle measures the deviation of a 2-plane from
holomorphicity.
A holomorphic 2-plane has $\Theta_{\rm K} = 0$, or 
$\Theta_{\rm K} = \pi$ (to see this 
consider a {\it holomorphic system of vectors ($J$-basis)} 
$A$, $\tilde{A}$ of this 2-plane,
$\tilde{\tilde{A}}= -A$)\footnote{For a holomorphic 
2-plane the action of the almost complex 
structure $J$ consist in a rotation within this 2-plane by the
(Euclidean) angle $\pi/2$ (= pseudo-angle, cf.\ sect.\ 4). 
In general, the almost complex 
structure $J$ maps an arbitrary 2-plane in $V_{\bf R}$ to a 2-plane
isoclinic to it (for an explanation of this term see sect.\ 5) 
with an angle which is equal to the K\"ahler angle
of the original 2-plane (\cite{roze3}, sect.\ 1, p.\ 65 
(p.\ 83 of the English translation),
\cite{gold}, chap.\ 1, sect.\ 1.3.3, p.\ 17).}, while 
an antiholomorphic 2-plane has $\Theta_{\rm K} = \frac{\pi}{2}$.
A 2-plane exhibiting some arbitrary K\"ahler angle $\Theta_{\rm K}$
is called $\Theta_{\rm K}${\it -holomorphic} \cite{tsuc}
or {\it general slant} \cite{chen1}, \cite{chen2}.
If one is just interested in the K\"ahler angle of a 2-plane
spanned by two vectors $A, B\in V_{\bf R}$
its defining equation (\ref{ZA6})
can be given a particularly simple shape if one chooses the vectors
$A$, $B$  such a way that they are orthogonal in the Euclidean sense 
($\Theta(A,B) = \frac{\pi}{2}$).\\

\noindent
{\bf 4.\ Relations between the different angles.}\ 
Now, applying unit vectors $a$, $b$ ($\vert c\vert =\vert
C\vert = \vert \tilde{C}\vert = 1$, $c = a,\ b$)
we can write eq.\ (\ref{ZA3}) as
\begin{eqnarray}
\label{ZA7}
\cos \Theta_{\rm c}(a,b)&=& 
\ \ (A,B)_{\bf R}\ \ +\  i\ \ (\tilde{A},B)_{\bf R}\nonumber\\[0.3cm]
&=&\cos\Theta (a,b) \ +\ i\ \cos\Theta_{\rm K}(a,b)\ \sin\Theta(a,b)
\ \ .
\end{eqnarray}
For the Hermitian angle $\Theta_{\rm H}$ defined by eq.\ (\ref{ZA5})
one immediately finds from eq.\ (\ref{ZA7}) the following relation 
to the Euclidean angle $\Theta$ and the K\"ahler angle $\Theta_{\rm K}$
(\cite{rizz4}, sect.\ 4, p.\ 97, eq.\ (15)).
\begin{eqnarray}
\label{ZA8}
\sin \Theta_{\rm H}(a,b)&=& 
\sin\Theta_{\rm K} (a,b) \ \sin\Theta (a,b)
\end{eqnarray}
The pseudo-angle $\varphi(a,b)$ can also be 
linked to the other angles defined above. Using eqs.\ (\ref{ZA4}), 
(\ref{ZA7}), one obtains the following
relations (also cf.\ \cite{gold}, chap.\ 2, sect.\ 2.2.2, 
pp.\ 36-38)\footnote{Eqs.\ (\ref{ZA9b}), (\ref{ZA9c}) have been written
in the most compact way ignoring the special cases to be treated
separately when the tangent or the cotangent functions become infinite.}.
\begin{eqnarray}
\label{ZA9a}
\cos \Theta (a,b)&=& 
\cos\Theta_{\rm H}(a,b)\ \cos\varphi (a,b)\\[0.3cm]
\label{ZA9b}
\sin \varphi(a,b)&=& 
\cot\Theta_{\rm K}(a,b)\ \tan\Theta_{\rm H} (a,b)\\[0.3cm]
\label{ZA9c}
\tan \varphi(a,b) &=& 
\cos\Theta_{\rm K}(a,b)\ \tan\Theta (a,b)
\end{eqnarray}
Specifically, one finds for any two vectors $a$, $b$ of a complex line
defining one and the same holomorphic 2-plane ($\Theta_{\rm K}(a,b) = 0$)
that the Hermitian angle vanishes
($\Theta_{\rm H}(a,b) = 0$, this follows from eq.\ (\ref{ZA8})) and 
that the pseudo-angle is equal to the Euclidean 
angle ($\varphi (a,b) = \Theta (a,b)$, this follows from eq.\ (\ref{ZA9c});
for $\Theta_{\rm K}(a,b) = \pi$ holds $\varphi (a,b) = - \Theta (a,b)$).
If the Hermitian angle $\Theta_{\rm H}(a,b)$ of two vectors $a$, $b$ is 
different from $\pi/2$ and if, in addition, these vectors are 
orthogonal in the Euclidean sense ($\Theta (a,b) = \frac{\pi}{2}$) 
their pseudo-angle $\varphi (a,b)$ is either equal to $\pi/2$
or $-\pi/2$ (cf.\ eqs.\ (\ref{ZA4}), (\ref{ZA7})). 
In the former case $\Theta_{\rm H} (a,b)
= \Theta_{\rm K} (a,b)$ applies while in the latter case
$\Theta_{\rm H} (a,b) = \pi - \Theta_{\rm K} (a,b)$
(cf.\ also \cite{roze3}, sect.\ 1, p.\ 65 
(p.\ 83 of the English translation)).
For any two vectors $a$, $b$ defining an antiholomorphic 2-plane
($\Theta_{\rm K}(a,b) = \frac{\pi}{2}$) from eq.\ (\ref{ZA7}) one
recognizes that the complex and the Euclidean angles 
coincide ($\Theta_{\rm c}(a,b) = \Theta (a,b)$) while for the
Hermitian angle holds $\Theta_{\rm H} (a,b) = \Theta (a,b)$,
or $\Theta_{\rm H} (a,b) = \pi - \Theta (a,b)$.
The pseudo-angle $\varphi (a,b)$ between vectors of an antiholomorphic 
2-plane is equal to $0$, or $\pm\pi$, respectively.\\

\noindent
{\bf 5.\ Further comments.}\ 
Before we continue the discussion with some further observations we need
to introduce the concept of isoclinic 2-planes. Two 2-planes 
$\cal A$, $\cal B$ in the real vector space 
$V_{\bf R}$ ($\simeq {\bf R}_{2n}$, $n\ge 2$) can intersect 
in various ways. In order to study their relation,
to each pair of lines ${\cal X}\subset {\cal A}$, 
${\cal Y}\subset {\cal B}$ 
the (Euclidean) angle they enclose can be calculated.
Once a line ${\cal X}\subset {\cal A}$ is fixed, for any arbitrary line 
${\cal Y}\subset {\cal B}$ the angle enclosed assumes values between some
$\alpha_0 \ge 0$  ($\alpha_0\le\frac{\pi}{2}$) 
and $\pi/2$. In general, $\alpha_0$ may lie 
between some minimal and some maximal value -- 
the so-called {\it stationary angles
(principal angles)} $\alpha_{\rm min}$,
$\alpha_{\rm max}$ -- which are characteristic for the 
geometry of the pair of 2-planes $\cal A$, $\cal B$. If $\alpha_{\rm min} =
\alpha_{\rm max}$, the 2-planes $\cal A$ and $\cal B$ are said to be
(mutually) {\it isoclinic} and the isocliny angle 
$\alpha_{\rm 0} = \alpha_{\rm min} = \alpha_{\rm max}$ of the two 2-planes 
can be determined by, say, calculating the (Euclidean) angle between 
some vector $X\in {\cal A}$ and its orthogonal projection onto $\cal B$
(for a comprehensive list of references on the subject of stationary 
angles and isoclinic subspaces see sect.\ II of \cite{scha}).\\

After these definitions we can now proceed with our discussion.
Any two holomorphic 2-planes
in $V_{\bf R}$ (corresponding to two complex lines in $V_{\bf C}$)
are isoclinic to each other (\cite{kwie}, \S 4, IX, p.\ 25,
\cite{maru1}, \S 3, p.\ 34,
\cite{wong2}, sect.\ 1-7, p.\ 51, theorem 1-7.4)\footnote{
This can easily be established by using holomorphic
systems of unit vectors for each 2-plane given by $a$ and $b$ 
($A$, $\tilde{A}$ and $B$, $\tilde{B}$, respectively) 
which are chosen such that $(A,B)_{\bf R} = 0$ holds. This way one
can calculate with little effort the orthogonal projection 
of a given unit vector of one 2-plane onto the other one and find
that the (Euclidean) angle between the vector and its projection 
is independent of the choice of the vector (the length (norm) of 
the projected vector is given by $\vert (\tilde{A},B)_{\bf R}\vert $).}. 
One can convince oneself that a vector of a holomorphic 2-plane 
in $V_{\bf R}$ and its (nonvanishing) orthogonal projection onto another 
(nonorthogonal) holomorphic 2-plane span an antiholomorphic 2-plane.
Therefore, the isocliny angle $\alpha_{\rm 0}$ between these 
two holomorphic 2-planes does not depend on the angle
concept applied for its definition as the complex, the Euclidean, and
the Hermitian angles agree for two lines of an antiholomorphic 2-plane.
As mentioned in sect.\ 4,
the Hermitian angle $\Theta_{\rm H}$ between two complex lines given
by two vectors $a$, $b$ which are orthogonal in the Euclidean sense
is determined by their K\"ahler angle through the equations
$\Theta_{\rm H} (a,b) = \Theta_{\rm K} (a,b)$  or
$\Theta_{\rm H} (a,b) = \pi - \Theta_{\rm K} (a,b)$. 
Therefore, if the two holomorphic 2-planes are given in terms of 
vectors $a$, $b$ which are orthogonal in the Euclidean sense
the isocliny angle between these 2-planes can be calculated 
by evaluating the K\"ahler angle of these vectors.\\

Finally, we would like to add some comments on the literature concerning
angles in vector spaces over fields. For complex vector spaces the 
relevant literature has been cited in course of the above discussion. 
The concept of the K\"ahler angle has also been extended to other,
in particular quaternionic, vector spaces \cite{roze5}, 
\cite{brun3}-\cite{brun10}, \cite{brun6}, \cite{brun7}, \cite{roze3},
\cite{roze7}, \cite{yukh}
(also see \cite{roze4}, sect.\ 6, p.\ 7 (p.\ 5 of the English 
translation)), \cite{roze2}, chap.\ III, \S 3.3, sect.\
3.3.6, pp.\ 182/183. Moreover, the concept of the K\"ahler angle
has also been generalized from 2-planes to higher dimensional
linear subspaces of real vector spaces which stand in 
correspondence to complex and other vector spaces \cite{rizz2},
\cite{rizz3}-\cite{ianu}. For further considerations of angles
in spaces over algebras see \cite{roze6}, chap.\ IV, p.\ 549,
\cite{bass}-\cite{herz}.\\

\small
{\bf Acknowledgements}. I am grateful to B.\ A.\ Rosenfeld for 
his comments on the first version of the paper and, in particular,
for pointing me to a paper by Goldman which meanwhile has appeared
in print as a book \cite{gold}.\\

\end{document}